\documentclass[11pt]{article}
\catcode`\@=11
\newtheorem{theorem}{Theorem}
\newtheorem{lemma}{Lemma}
\newtheorem{proposition}{Proposition}
\newtheorem{definition}{Definition}
\newtheorem{example}{Example}
\newtheorem{remark}{Remark}

\newtheorem{corollary}{Corollary}
\def\demo{\noindent{\bf Proof .-}}
\def\section{\@startsection {section}{1}{\z@}{-3.5ex plus -1ex
minus-.2ex}{2.3ex plus .2ex}{\normalsize\bf}}

\def\bn{\hbox{\it I\hskip -2pt N}}
\def\bz{\hbox{\it Z\hskip -4pt Z}}

\def\bp{\hbox{\it P\hskip -5.5pt P}}

\begin{document}
\begin{center}
{\Large\bf \textsc{On binomial set-theoretic complete intersections in characteristic $p$}}
\end{center}
\vskip.5truecm
\begin{center}
{Margherita Barile\footnote{Partially supported by the Italian Ministry of Education, University and Research.}\\ Dipartimento di Matematica, Universit\`{a} di Bari,Via E. Orabona 4,\\70125 Bari, Italy}
\end{center}
\vskip1truecm
\noindent
{\small {\bf Abstract} Using aritmethic conditions on affine semigroups we prove that for a simplicial toric variety of codimension 2 the property of being a set-theoretic complete intersection on binomials in characteristic $p$ holds either for all primes $p$, or for no prime $p$, or for exactly one prime $p$. \vskip0.5truecm
\noindent
}

\section*{Introduction} It has been recently discovered that the minimum number of equations that are needed to define an algebraic variety, in general, depends on the characteristic of the ground field. The first examples to be found were the determinantal and Pfaffian varieties considered in \cite{B0}. Later, in \cite{B1}, \cite{B3}, \cite{B4} and \cite{BL}, we presented infinite classes of toric varieties which are set-theoretic complete intersections (i.e., are defined by the vanishing of as many polynomials as their codimension) only in one positive characteristic. In all these cases, the polynomials were, in fact, binomials. This raised the question whether there exists a toric variety $V$ which is defined by codim\,$V$ binomial equations in more than one, but not in all positive characteristics. In the present paper we show that the answer is negative if $V$ is simplicial and codim\,$V=2$. We give a complete classification of the possible cases in terms of the arithmetic properties of the affine semigroups attached to these varieties. This is an exhaustive study, which includes and completes all previously known results.    
\section{Preliminaries}
Let $K$ be an algebraically closed field, and let $T$ be a finite subset of $\bz^n$ $(n\geq2)$ having cardinality $m>0$, say $T=\{{\bf t}_1, \dots, {\bf t}_m\}$, where ${\bf t}_i=(t_{i1}, \dots, t_{in})$ for all $i=1,\dots, m$. We will throughout assume that the set $T$ is not free over $\bz$. Consider the homomorphism of polynomial rings
$$\phi_T:K[z_1,\dots, z_m]\to K[u_1,\dots, u_n]$$
\noindent such that
$$\phi_T(z_i)=u_1^{t_{i1}}\cdots u_n^{t_{in}},\qquad\qquad\mbox{for all }i=1,\dots, m.$$
\noindent
The ideal $I_T=\ker\phi_T$ is generated by all binomials 
\begin{equation}\label{binomials}B^{\alpha_1^+\cdots\alpha_m^+}_{\alpha_1^-\cdots\alpha_m^-}=
z_1^{\alpha_1^+}\cdots z_m^{\alpha_m^+}-
z_1^{\alpha_1^-}\cdots z_m^{\alpha_m^-}\end{equation}
\noindent
where $\alpha_i^+,\alpha_i^-$ are nonnegative integers (not all zero) such that
\begin{equation}\label{relations}\alpha_1^+{\bf t}_1+\cdots+\alpha_m^+{\bf t}_m=\alpha_1^-{\bf t}_1+\cdots+\alpha_m^-{\bf t}_m,\end{equation}
\noindent
and in fact there is a one-to-one correspondence between the binomials (\ref{binomials}) in $I_T$ and the semigroup relations (\ref{relations}) in $\bn T$.
The ideal $I_T$ is the defining ideal of the variety $V_T\subset K^m$ parametrized by
$$z_1=u_1^{t_{11}}\cdots u_n^{t_{1n}},\qquad \dots,\qquad  z_m=u_1^{t_{m1}}\cdots u_n^{t_{mn}},$$
\noindent
which is an affine toric variety.  If all sums $\sum_{j=1}^nt_{ij}$, for $i=1,\dots, m$, are equal, then this parametrization is called {\it homogeneous} and it defines a projective toric variety in ${\bp}^{m-1}$. The dimension of $V_T$ is equal to the rank of the abelian group $\bz T$ if $V$ is affine, it is one less if $V$ is projective. \newline
If, $n<m$, and,  up to a change of parameters, for all $i=1,\dots, n$ we have ${\bf t}_i=c{\bf e}_i$ with $c$ a positive integer and ${\bf e}_i$ the $i$th element of the canonical basis of ${\bz}^n$, then $V$ is called {\it simplicial}, and we have $\dim V=n$ if $V$ is affine, $\dim V=n-1$ if $V$ is projective. In this case, whenever $t_{ij}\ne0$ for all $i=n+1, \dots, m$ and all $j=1,\dots, n$, the parametrization of $V$ is called {\it full}. A simplicial toric variety of dimension 1 is called a {\it monomial curve}: we have an affine monomial curve for $n=1$ and a projective monomial curve if the parametrization is homogeneous and $n=2$; in both cases the parametrization is obviously full.  
Set codim\,$V_T=r$. Then we say that $V_T$ is a {\it (binomial) set-theoretic complete intersection} if there are (binomials) $F_1,\dots, F_r\in K[z_1,\dots, z_m]$ such that 
$$I_T= \sqrt{(F_1,\dots, F_r)},$$
\noindent
which, by Hilbert's Nullstellensatz, is equivalent to stating that $V_T$ is defined by the following system of $r$ (binomial) equations:
$$F_1=\cdots=F_r=0.$$
\noindent
According to the above description, for the variety $V_T$ the property of being a complete intersection (i.e., the equality between the minimal number of generators of $I_T$ and codim\,$V_T$) only depends on $T$, and is thus independent of the ground field $K$. The (binomial) set-theoretic complete intersection property, however, depends, in general, on char\,$K$. Barile, Morales and Thoma proved the following result.
\begin{proposition}{\rm (\cite{BMT2}, Theorem 4)}\label{Barile} If char\,$K=0$, then $V_T$ is a binomial set-theoretic complete intersection if and only if it is a complete intersection. 
\end{proposition}
\begin{remark}\label{remark1}  {\rm There are many examples of affine or projective monomial curves which are not complete intersections, but are set-theoretic complete intersections in characteristic zero: some of them can be found in \cite{Br1}, \cite{Br2}, \cite{E}, \cite{K}, \cite{RV1}, \cite{RV2}, \cite{S},  \cite{T} and \cite{V}.  This fact, together with Proposition \ref{Barile}, shows that being a set-theoretic complete intersection is, in general, a strictly weaker condition than being a \textsl{binomial} set-theoretic complete intersection. We do not know, however, whether this is true in positive characteristics.}
\end{remark}
The  sets $T$ for which the variety $V_T$ is a complete intersection have been completely characterized by means of the following notions, introduced by Rosales \cite{R}.
\begin{definition}\label{gluing} {\rm
Let $T_1$ and $T_2$ be nonempty
subsets of $T$ such that $T     =   T_1\cup T_2$ and $T_1\cap T_2    =
\emptyset$.
Then $T$ is called a {\it gluing} of $T_1$ and
$T_2$ if  there is  a
nonzero element ${\bf w}\in {\bn}^n$ such that  ${\bz} T_1\cap{\bz} T_2    =  {\bz} {\bf w}$
and 
${\bf w}\in {\bn } T_1 \cap {\bn } T_2$.}
\end{definition}
\begin{definition}\label{ci}{\rm
The set $T$ is called a {\it complete intersection} if $T$
is the
gluing of
$T_1$ and $T_2$, where each  of the subsets $T_1, 
T_2$
is a complete intersection or generates a free abelian semigroup.}
\end{definition}
The next result is due to Fischer, Morris and Shapiro \cite{FMS}.
\begin{proposition}\label{Fischer} The variety $V_T$ is a complete intersection if and only if $T$ is the gluing of two subsets $T_1$ and $T_2$ which are complete intersections.
\end{proposition}
Proposition \ref{Barile} is false in positive characteristic: here a necessary and sufficient condition for being a binomial set-theoretic complete intersection  can be formulated in terms of an arithmetic property of the set $T$ and of the affine semigroup $\bn T$, which is derived from the ones given in Definitions \ref{gluing} and \ref{ci} and was introduced in 
 \cite{BMT2}, pp.~1894--1895, by means of the following two definitions.
 \begin{definition}\label{pgluing}{\rm
Let $p$ be a prime number and let $T_1$ and $T_2$ be nonempty
subsets of $T$ such that $T     =   T_1\cup T_2$ and $T_1\cap T_2    =
\emptyset$.
Then $T$ is called a {\it p-gluing} of $T_1$ and
$T_2$ if  there are a nonnegative integer $k$ and  a
nonzero element ${\bf w}\in {\bz}^n$ such that  ${\bz} T_1\cap{\bz} T_2    =  {\bz} {\bf w}$
and 
$  p^k{\bf w}\in {\bn } T_1 \cap {\bn } T_2$.}
\end{definition}
\begin{definition}{\rm
The set $T$ is called {\it completely $p$-glued} if $T$
is the
$p$-gluing of
$T_1$ and $T_2$, where each  of the subsets $T_1, 
T_2$
is completely $p$-glued or generates a free abelian semigroup.}
\end{definition}
\begin{remark}\label{remark0}{\rm The comparison of the above definitions yields the following.
\begin{list}{}{}
\item{a)} If $T$ is the gluing of two nonempty subsets $T_1$ and $T_2$, then it is the $p$-gluing of $T_1$ and $T_2$ for all primes $p$.
\item{b)} If $T$ is a complete intersection, then it is completely $p$-glued for all primes $p$.
\end{list}
}
\end{remark}
\begin{proposition}\label{p}{\rm (\cite{BMT2}, Theorem 5)} If char\,$K=p$, then $V_T$ is a binomial set-theoretic complete intersection if and only if $T$ is completely $p$-glued. 
\end{proposition}
The papers \cite{B1}, \cite{B3}, \cite{B4}, \cite{BL} present classes of toric varieties that fulfil this property  in exactly one positive characteristic $p$, whereas those described in \cite{B2}, and some of those contained in \cite{B4}, do not fulfil it for any $p$. In all other known cases the property holds for all $p$. We recall two general results in this direction.
\begin{proposition}\label{full}{\rm (\cite{BMT1}, Theorem 1)} Every simplicial toric variety with a full par\-am\-etriz\-ation is a binomial set-theoretic intersection in all positive characteristics.  
\end{proposition}
This proposition includes the following one as a special case.
\begin{proposition}\label{moh}{\rm (\cite{M}, Corollary)} Every projective monomial curve is a binomial set-theoretic complete intersection in all positive characteristics.  
\end{proposition}
\noindent
So far no toric variety is known which is a binomial set-theoretic complete intersection in more than one, but not in all positive characteristics. The aim of this paper is to show that no such examples exist among the simplicial toric varieties of codimension 2. In Sections 3--5 we will give a complete classification of $p$-gluings for affine  simplicial toric varieties of codimension 2. In Section 4, in particular, we will introduce an arithmetic invariant which allows us to distinguish between the possible cases.  \newline
In Section 2 we will preliminarily clarify the role of  Propositions \ref{full} and \ref{moh} in the general arithmetical framework of completely $p$-glued sets. We will also prove the following additional result.
\begin{proposition}\label{additional} Every affine simplicial toric variety of dimension at most 2 is a binomial set-theoretic complete intersection in all positive characteristics.
\end{proposition} 
\section{Some known results revisited}
In the sequel, given any ${\bf v}\in\bz^n$, we will denote by $v_i$ its $i$-th component, for $1\leq i\leq n$.\newline
We will consider the varieties $V_T$, where $T=\{{\bf a}, {\bf b}, c{\bf e}_1, \dots, c{\bf e}_n\}\subset \bz^n$, and $c$ is any positive integer. This class consists of all affine (or projective) simplicial  toric varieties of codimension 2 in $K^{n+2}$ (or in $\bp^{n+1}$). It is not restrictive to assume that, for all $i=1,\dots, n$, we have $a_i\ne 0$ or $b_i\ne 0$. We introduce the following two partitions of $T$:
$$\tilde T_1=\{{\bf a}\},\qquad \tilde T_2=\{{\bf b}, c{\bf e}_1, \dots, c{\bf e}_n\},$$
$$\tilde{\tilde T_1}=\{{\bf b}\},\qquad \tilde{\tilde T_2}=\{{\bf a}, c{\bf e}_1, \dots, c{\bf e}_n\}.$$
\noindent
For all ${\bf v}\in {\bz}^n$ we will set supp\,${\bf v}=\{i\vert v_i\ne 0\}.$
The content of the next two propositions is not totally new, since they partly consist of results collected from other sources. Nevertheless, we find it useful and interesting to re-propose these statements here in a unified form, which emphasizes the role of $p$-gluings. The proof of the next proposition under the assumption (i) is hinted at in \cite{BMT2}, Example 1. 
\begin{proposition}\label{proposition1} Suppose that either
\begin{list}{}{}
\item{(i)}  {\rm supp}\,${\bf b}\subset\,{\rm supp}\,{\bf a}$,
or
\item{(ii)}  {\rm supp}\,${\bf b}=\{1,\dots, n\}\setminus\,{\rm supp}\,{\bf a}$.
\end{list}
\noindent
Then $T$ is the $p$-gluing of $\tilde T_1$ and $\tilde T_2$ for all primes $p$.\newline
(If ${\bf a}$ and ${\bf b}$ are exchanged, then $T$ is the $p$-gluing of $\tilde{\tilde T_1}$ and $\tilde{\tilde T_2}$ for all primes $p$.)
\end{proposition}
\demo Suppose that (i) holds. Since 
$$c{\bf a}=\sum_{i=1}^n a_ic{\bf e}_i\in{\bz}\tilde T_2,$$
\noindent
the number
\begin{equation}\label{alpha}\alpha=\gcd\{\lambda\in\bz\vert \lambda{\bf a}\in\bz\tilde T_2\}\end{equation}
\noindent
is a well-defined positive integer. Moreover 
$$\bz\tilde T_1\cap\bz\tilde T_2=\bz\alpha{\bf a}.$$
\noindent
 We have that, for some integers $\beta, \gamma_1, \dots, \gamma_n$,
$$\alpha{\bf a}=\beta{\bf b}+\sum_{i=1}^n\gamma_ic{\bf e}_i,$$
\noindent
i.e., for all $i=1,\dots, n$ we have the following equality, where either side is a natural number:
\begin{equation}\label{0}
\alpha a_i=\beta b_i+\gamma_ic.\end{equation}
\noindent
We have to show that under the given assumption, for all primes $p$ there is a nonnegative integer $k$ such that 
\begin{equation}\label{ast} p^k\alpha{\bf a}\in\bn\tilde T_2.\end{equation}
\noindent 
Suppose that {\rm supp}\,${\bf b}\subset\,{\rm supp}\,{\bf a}$. Then, for all $i\in\,{\rm supp}\,{\bf b}$, in (\ref{0}) we have  $\alpha a_i >0$, whence $-\frac{\gamma_i}{b_i}<\frac{\beta}c$. Set 
\begin{equation}\label{1}\gamma=\max_{i\in\,{\rm supp}\,{\bf b}}-\frac{\gamma_i}{b_i}.\end{equation}
\noindent Then $\gamma<\frac{\beta}c$. 
Let $p$ be any prime number. We can choose a sufficiently large positive integer $k$ so as to have 
\begin{equation}\label{2}
p^k\gamma<d<p^k\frac{\beta}c
\end{equation}
for some integer $d$. Set
$$\beta'=p^k\beta-dc,$$
\noindent
and, for all $i=1,\dots, n$, 
$$\gamma_i'=p^k\gamma_i+db_i.$$
\noindent
Then (\ref{2}) implies that $\beta'>0$ and, moreover, from (\ref{1}) and (\ref{2}), for all $i\in\,{\rm supp}\,{\bf b}$, we get 
$$p^k\frac{\gamma_i}{b_i}\geq -p^k\gamma>-d,$$
\noindent
which implies that $\gamma_i'>0.$ If $i\not\in\,{\rm supp}\,{\bf b}$, then from (\ref{0}) we deduce that 
$\alpha a_i=\gamma_ic,$
which gives us $\gamma_i>0$, so that $\gamma_i'=p^k\gamma_i>0$. Hence 
$$p^k\alpha {\bf a}= \beta'{\bf b}+\sum_{i=1}^n\gamma_i'c{\bf e}_i,$$
\noindent
so that (\ref{ast}) is true. \newline
Now suppose that (ii) holds. For all $i\in\,{\rm supp}\,{\bf a}$ we have  $a_i>0$ and $b_i=0$. Thus, in view of (\ref{0}) we also have $\alpha a_i=\gamma_i c,$
where $\gamma_i>0$. Therefore 
$$\alpha {\bf a}=\sum_{i\in\,{\rm supp}\,{\bf a}}\gamma_ic{\bf e}_i\in\bn \tilde T_2,$$
so that (\ref{ast}) is fulfilled for $k=0$. This completes the proof.
\par\smallskip\noindent
\begin{proposition}\label{proposition2} a) If $T$ is the gluing of $\tilde T_1$ and $\tilde T_2$ (or of $\tilde{\tilde T_1}$ and $\tilde{\tilde T_2}$), then $T$ is a complete intersection.\newline
b)  Let $p$ be a prime. If $T$ is the $p$-gluing of $\tilde T_1$ and $\tilde T_2$ (or of $\tilde{\tilde T_1}$ and $\tilde{\tilde T_2}$), then $T$ is completely $p$-glued. 
\end{proposition}
\demo  We show the claim for $\tilde T_1$ and $\tilde T_2$. In view of Remark \ref{remark0} a), it suffices to show that $\tilde T_2$ is the gluing of $\tilde T_{21}=\{{\bf b}\}$ and $\tilde T_{22}=\{c{\bf e}_1,\dots, c{\bf e}_n\}$: in fact, $\bn \tilde T_{21}$ and $\bn \tilde T_{22}$ are  both  free abelian semigroups. Now, if $\beta=c/\gcd(c,b_1,\dots, b_n)$, then $\bz\beta{\bf b}=\bz\tilde T_{21}\cap\bz\tilde T_{22}$  and $\beta{\bf b}\in\bn\tilde T_{21}\cap\bn\tilde T_{22}$. The claim for $\tilde{\tilde T_1}$ and $\tilde{\tilde T_2}$ can be shown similarly. This completes the proof. 
\par\smallskip\noindent
\begin{corollary}\label{corollary1} If $n\leq 2$, then, for all primes $p$, $T$ is the $p$-gluing of $\tilde T_1$ and $\tilde T_2$, or of $\tilde{\tilde T_1}$ and $\tilde{\tilde T_2}$. In particular, it is completely $p$-glued for all primes $p$.  \end{corollary}
\demo The first part of the claim follows immediately from Proposition \ref{proposition1}, since for $n=2$ one of (i) and (ii) must hold. The second part is a consequence of Proposition \ref{proposition2} b).
\par\smallskip\noindent
\begin{remark}{\rm Proposition \ref{full}, Proposition \ref{moh} and Proposition \ref{additional} follow from Corollary \ref{corollary1} and Proposition \ref{p}. 
}
\end{remark}
\section{The main theorem on $p$-gluings} In this section we prove the crucial technical result. We consider the set $T$ introduced above and assume that it does not fulfil the assumption (i) of Proposition \ref{proposition1} with respect to any order of ${\bf a}$ and ${\bf b}$. We show that the property of being $p$-glued then reduces to the property of being the $p$-gluing of $\tilde T_1$ and  $\tilde T_2$ (or, equivalently, of $\tilde{\tilde{ T_1}}$ and  
$\tilde{\tilde{ T_2}}$).
\begin{theorem}\label{main} Let $p$ be a prime. Suppose that none of supp\,${\bf a}$ and supp\,${\bf b}$ is contained in the other. If $T$ is the $p$-gluing of two nonempty subsets $T_1$ and $T_2$, then  $T$ is the $p$-gluing of $\tilde T_1$ and  $\tilde T_2$, and of $\tilde{\tilde{ T_1}}$ and  
$\tilde{\tilde{ T_2}}$.
\end{theorem}
\demo  As a trivial consequence of the first assumption, we have  $n\geq 2$, and we may assume that, up to a change of indices, $a_1\ne0$, $b_1=0$, $a_2=0$, $b_2\ne0$. By the second assumption there is ${\bf w}\in\bz^n$, ${\bf w}\ne0$, such that 
\begin{equation}\label{8'} \bz{\bf w}=\bz T_1\cap\bz T_2\qquad\mbox{and}\qquad p^k{\bf w}\in\bn T_1\cap\bn T_2\end{equation}
\noindent
for some nonnegative integer $k$. Let 
$$T_1=\{{\bf t}_{11},\dots, {\bf t}_{1r}\},\qquad \mbox{and}\qquad  T_2=\{{\bf t}_{21},\dots, {\bf t}_{2s}\}.$$
\noindent Then we have 
\begin{equation}\label{3} {\bf w}=\sum_{i=1}^r\alpha_i{\bf t}_{1i}=\sum_{i=1}^s\beta_i{\bf t}_{2i}\end{equation}
\noindent
for some 
$\alpha_1,\dots, \alpha_r,\beta_1, \dots, \beta_s\in\bz,$ and
\begin{equation}\label{4} p^k{\bf w}=\sum_{i=1}^r\tilde\alpha_i{\bf t}_{1i}=\sum_{i=1}^s\tilde\beta_i{\bf t}_{2i}\end{equation}
\noindent
for some 
$\tilde\alpha_1,\dots, \tilde\alpha_r,\tilde\beta_1, \dots, \tilde\beta_s\in\bn$. Since $T_1\cup T_2=T$, up to exchanging $T_1$ and $T_2$ we may assume that ${\bf a}\in T_1$, say ${\bf t}_{11}={\bf a}$, so that $\alpha_1$ is the coefficient of ${\bf a}$ in (\ref{3}) and $\tilde\alpha_1$ is the coefficient of ${\bf a}$ in (\ref{4}). We will also adopt the following notation:
$$I_1=\{i\vert c{\bf e}_i\in T_1\}\qquad\mbox{and}\qquad I_2=\{i\vert c{\bf e}_i\in T_2\}.$$
\noindent
The theorem will be proven in several steps, through a sequence of intermediate claims. 
\par\smallskip\noindent
\underline{Claim 1} a) One of ${\bf a}$ and ${\bf b}$ appears in (\ref{4}) with a nonzero coefficient.\newline
b)  If $\tilde\alpha_1\ne0$, then $\alpha_1\ne 0$ and $1\in I_2$. 
\newline Proof of Claim 1. If ${\bf a}$ and ${\bf b}$ both appeared in (\ref{4}) with a zero coefficient, then we would have ${\bf w}=0$, because the set $\{c{\bf e}_1, \dots, c{\bf e}_n\}$ is free, against our assumption on ${\bf w}$. This proves a). We now prove b). Suppose that $\tilde\alpha_1\ne 0$. Then $1\in I_2$, otherwise we could assume that ${\bf t}_{12}=c{\bf e}_1$ and we would deduce that
$$p^kw_1=\tilde\alpha_1a_1+\tilde\alpha_2c=0,$$
\noindent
because $b_1=0$. Since all elements in the above equality are nonnegative integers, it would follow that $a_1=0$, against our assumption. Now suppose for a contradiction that $\alpha_1=0$. Then, since $1\in I_2$ and $b_1=0$, from (\ref{3}) we deduce that $w_1=0$. But the element 
\begin{equation}\label{u}{\bf u}=c{\bf a}-\sum_{i\in I_1}a_ic{\bf e}_i=\sum_{i\in I_2}a_ic{\bf e}_i\end{equation}
 \noindent
also belongs to $\bz T_1\cap\bz T_2$. However, since $u_1=a_1c\ne0$, we conclude that
${\bf u}$ cannot belong to $\bz{\bf w}$, which contradicts (\ref{8'}). This shows that $\alpha_1\ne0$ and completes the proof of Claim 1.
\par\smallskip\noindent
\underline{Claim 2} If $1\in I_2$, then 
$$\bz\alpha_1{\bf a}=\bz\tilde T_1	\cap\bz\tilde T_2,$$
\noindent
and $\tilde\alpha_1=p^k\alpha_1$. 
\newline
Proof of Claim 2. In view of (\ref{3}),
$$\alpha_1{\bf a}\in\bz\tilde T_1\cap\bz\tilde T_2.$$
\noindent
This intersection is a subgroup of $\bz{\bf a}$, hence it is $\bz\gamma{\bf a}$, for some integer $\gamma$ having the same sign as $\alpha_1$. Since $\alpha_1{\bf a}\in\bz\gamma{\bf a}$, we also have that $\gamma$ divides $\alpha_1$. Since $\gamma{\bf a}\in\bz\tilde T_2$, the element $\gamma{\bf a}=\gamma{\bf t}_{11}$ can be expressed as a linear combination, with integer coefficients, of  ${\bf b}, {\bf e}_1,\dots, {\bf e}_n$, i.e., of ${\bf t}_{12},\dots, {\bf t}_{1r}, {\bf t}_{21}, \dots, {\bf t}_{2s}$. We thus have an equality
\begin{equation}\label{7} \sum_{i=1}^r\alpha_i'{\bf t}_{1i}=\sum_{i=1}^s\beta_i'{\bf t}_{2i}\end{equation}
\noindent
for $\alpha_1'=\gamma$ and some $\alpha_2',\dots, \alpha_r',\beta_1',\dots, \beta_s'\in\bz$. The element on either side of (\ref{7}) belongs to $\bz T_1\cap\bz T_2$, hence it is equal to $\lambda{\bf w}$ for some $\lambda\in\bz$. Thus (\ref{7}) implies that
\begin{equation}\label{8} \lambda w_1=\gamma a_1,\end{equation}
because $1\in I_2$ and $b_1=0$. On the other hand, (\ref{3}) yields
\begin{equation}\label{9} w_1=\alpha_1a_1,\end{equation}
\noindent
but then, replacing (\ref{9}) in (\ref{8}), we get $\lambda\alpha_1=\gamma$, so that $\alpha_1$ divides $\gamma$. We thus conclude that $\alpha_1$ and $\gamma$ are associate integers of the same sign, i.e., $\alpha_1=\gamma$. Furthermore, equality (\ref{4}) gives us $p^kw_1=\tilde\alpha_1a_1,$ from which, in view of (\ref{9}), we deduce that $\tilde\alpha_1=p^k\alpha_1$, as required. This completes the proof of Claim 2. 
\par\smallskip\noindent
\underline{Claim 3} If $\tilde\alpha_1\ne0$, then ${\bf b}$ appears in (\ref{4}) with coefficient zero.
\newline
Proof of Claim 3. Suppose for a contradiction that, under the given hypothesis, ${\bf b}$ appears in (\ref{4}) with a nonzero coefficient. There is $j\in\{1,2\}$ such that ${\bf b}\in T_j$. Therefore, replacing ${\bf a}$ with ${\bf b}$, ${\bf e}_1$ with ${\bf e}_2$, and $T_1$ with $T_j$, Claim 1 b) allows us to conclude that $c{\bf e}_2\not\in T_j$. Then
$${\bf v}=c{\bf b}-\sum_{i\in I_j}b_ic{\bf e}_i=\sum_{i\not\in I_j}b_ic{\bf e}_i\in{\bz}T_1\cap\bz T_2,$$
\noindent
and $v_1=0$, because $b_1=0$, whereas  $v_2=b_2c\ne0$.  But then ${\bf v}$ and the element ${\bf u}$ from (\ref{u}) cannot both belong to $\bz{\bf w}$, which contradicts (\ref{8'}). This proves Claim 3.
\par\smallskip\noindent
\underline{Claim 4} If $\tilde\alpha_1\ne0$, then $T$ is the $p$-gluing of $\tilde T_1$ and $\tilde T_2$.
\newline
Proof of Claim 4. From (\ref{4}) we get, in view of Claim 1 b), Claim 2 and Claim 3,
$$\tilde\alpha_1{\bf a}=p^k\alpha_1{\bf a}\in\bz\{c{\bf e}_1, \dots, c{\bf e}_n\},$$
\noindent
which implies that
\begin{equation}\label{10} p^k\vert\alpha_1\vert{\bf a}\in\bn\{c{\bf e}_1,\dots, c{\bf e}_n\}\subset\bn\tilde T_2.\end{equation}
\noindent
On the other hand, by Claim 2, 
\begin{equation}\label{11} \bz\vert\alpha_1\vert{\bf a}=\bz\tilde T_1\cap\bz\tilde T_2,\end{equation}
\noindent
which, together with (\ref{10}), implies that $T$ is the $p$-gluing of $\tilde T_1$ and $\tilde T_2$.   This proves Claim 4.
\par\smallskip\noindent
\underline{Claim 5} $T$ is the $p$-gluing of $\tilde{\tilde T_1}$ and $\tilde{\tilde T_2}$.
\newline
Proof of Claim 5. According to Claim 1 a), we have that either ${\bf a}$ or ${\bf b}$ appears in (\ref{4}) with a nonzero coefficient. In the latter case we have that $T$ is the $p$-gluing of $\tilde{\tilde T_1}$ and $\tilde{\tilde T_2}$: this follows from Claim 4, after exchanging ${\bf a}$ and ${\bf b}$.   So assume that the former case occurs, i.e., that $\tilde\alpha_1\ne0$. Then, by Claim 1 b), $1\in I_2$. Let $\beta'\geq 0$ be such that 
\begin{equation}\label{13}\bz\beta'{\bf b}=\bz\tilde{\tilde T_1}\cap\bz\tilde{\tilde T_2}.\end{equation}
\noindent
Then, for some integers $\alpha', \gamma_1', \dots, \gamma_n'$, 
\begin{equation}\label{14} \beta'{\bf b}=\alpha'{\bf a}+\sum_{i=1}^n\gamma_i'c{\bf e}_i,
\end{equation}
so that, $\alpha'{\bf a}\in\bz\tilde T_1\cap\bz\tilde T_2$ and therefore, in view of Claim 2, $\alpha_1$ divides $\alpha'$. Hence (\ref{14}) implies that
$$\beta'{\bf b}\in\bz\{\alpha_1{\bf a}, c{\bf e}_1,\dots, c{\bf e}_n\},$$
\noindent
whence
$$p^k\beta'{\bf b}\in\bz\{p^k\alpha_1{\bf a}, c{\bf e}_1,\dots, c{\bf e}_n\}\subset\bz\{c{\bf e}_1, \dots, c{\bf e}_n\},$$
\noindent
where the inclusion follows from (\ref{10}). Therefore $p^k\beta'{\bf b}\in\bn\{c{\bf e}_1,\dots, c{\bf e}_n\}\subset\bn\tilde{\tilde T_2}.$ This, together with (\ref{13}), shows that $T$ is the $p$-gluing of $\tilde{\tilde T_1}$ and $\tilde{\tilde T_2}$. This completes the proof of Claim 5.\par\smallskip\noindent
We can now conclude the proof of the theorem. By Claim 5, $T$ is the $p$-gluing of $\tilde{\tilde T_1}$ and $\tilde{\tilde T_2}$. The argumentations developed above, with ${\bf a}$ and ${\bf b}$ exchanged, allow us to conclude that $T$ is the $p$-gluing of $\tilde T_1$ and $\tilde T_2$ as well. This completes the proof of the theorem. 
\par\smallskip\noindent 
\begin{lemma}\label{lemma1} If $T$ is the $p$-gluing of two nonempty subsets $T_1$ and $T_2$ for two different primes $p$, then it is the $p$-gluing of $T_1$ and $T_2$ for all primes $p$.
\end{lemma}
\demo If $T$ is the $p$-gluing and the $p'$-gluing of $T_1$ and $T_2$, where $p$ and $p'$ are distinct primes, then there is ${\bf w}\in\bz^n$ such that $\bz{\bf w}=\bz T_1\cap\bz T_2$ and
\begin{equation}\label{12} p^k{\bf w}, p'^{k'}{\bf w}\in\bn T_1\cap\bn T_2,\end{equation}
\noindent
for some nonnegative integers $k, k'$. Let $q$ be a prime. Then there are $a,b,\ell\in\bn$ such that
$$q^{\ell}=ap^k+bp'^{k'},$$
\noindent
so that, by (\ref{12}), 
$$q^{\ell}{\bf w}\in\bn T_1\cap\bn T_2.$$
\noindent
This proves that $T$ is the $q$-gluing of $T_1$ and $T_2$, which completes the proof. 
\par\smallskip\noindent
Proposition \ref{proposition1}, Theorem \ref{main} and Lemma \ref{lemma1} imply the following result.
\begin{corollary}\label{corollary6} If $T$ is completely $p$-glued for two distinct primes $p$, then it is completely $p$-glued for all primes $p$. 
\end{corollary}
The converse of statement b) of Remark \ref{remark0} is not true: Propositions \ref{full} and \ref{moh} provide infinitely many examples of sets $T$ which are completely $p$-glued for all primes $p$, but are not complete intersections. We, however, have the following result.
\begin{corollary}\label{corollary6bis} Suppose that neither of supp\,${\bf a}$ and supp\,${\bf b}$ is contained in the other. Then $T$ is completely $p$-glued for all primes $p$ if and only if $T$ is a complete intersection. 
\end{corollary}
\demo   We just need to prove the {\it only if} part.  Suppose that $T$ is completely $p$-glued and $q$-glued, where $p$ and $q$ are distinct primes.  By virtue of Theorem \ref{main},  $T$ is the $p$-gluing and the $q$-gluing of $\tilde T_1$ and $\tilde T_2$, so that, for some nonzero $\alpha\in\bn$, we have
\begin{equation}\label{I}\bz\alpha{\bf a}=\bz\tilde T_1\cap\bz\tilde T_2,\end{equation}
\noindent
 and $p^k\alpha{\bf a}\in\bn\tilde T_2$, $q^k\alpha{\bf a}\in\bn\tilde T_2$,  for some nonnegative integer $k$. Since  supp\,${\bf b}\not\subset\,$\,supp\,${\bf a}$, this implies that $p^k\alpha{\bf a}, q^k\alpha{\bf a}\in\bn\{c{\bf e}_1,\dots, c{\bf e}_n\}$. Let $g=c/\gcd(c, a_1, \dots, a_n)$. Then ${\bn} g{\bf a}={\bn}{\bf a}\cap{\bn}\{c{\bf e}_1,\dots, c{\bf e}_n\}$. Thus $g$ divides both $p^k\alpha$ and $q^k\alpha$, so that $g$ divides $\alpha$. But then  $\alpha{\bf a}\in{\bn}\{c{\bf e}_1,\dots, c{\bf e}_n\}$, whence 
\begin{equation}\label{II}\alpha{\bf a}\in{\bn}{\tilde T_1}\cap{\bn}{\tilde T_2}.\end{equation}
\noindent
  Relations (\ref{I}) and (\ref{II}) show that $T$ is the gluing of $\tilde T_1$ and $\tilde T_2$. In view of Proposition \ref{proposition2} a), it follows that $T$ is a complete intersection. This completes the proof.
\par\smallskip\noindent
 Proposition \ref{proposition1}, Theorem  \ref{main} and Corollary \ref{corollary6} imply the following. 
\begin{corollary}\label{corollary_added} The set $T$ is completely $p$-glued
\begin{list}{}{}
\item{(a)} for exactly one prime $p$, or
\item{(b)} for  every prime $p$, or
\item{(c)} for no prime $p$.
\end{list}
\noindent
Moreover, if $T$ is completely $p$-glued, then it is the $p$-gluing of $\tilde T_1$ and $\tilde T_2$ (or of $\tilde{\tilde T_1}$ and $\tilde{\tilde T_2}$).
\end{corollary}
\section{A computational point of view}
In this section we provide a complete arithmetic characterization of the sets $T$ fulfilling each of the cases listed in Corollary \ref{corollary_added} under the assumption that neither of supp\,${\bf a}$ and supp\,${\bf b}$ is contained in the other. Our first step is to derive a closed formula for the number $\alpha$ in (\ref{alpha}). Up to dividing the elements of $T$ by their greatest common divisor $\delta$ (which corresponds to a change of parameters for $V_T$, in which, for all $i=1\,\dots, n$, $u_i^c$, $u_i^{a_i}$ and $u_i^{b_i}$ are replaced by $u_i^{c/\delta}$, $u_i^{a_i/\delta}$ and $u_i^{b_i/\delta}$, respectively), we may assume that 
\begin{equation}\label{gcd} \gcd(c, a_1,\dots, a_n,b_1,\dots, b_n)=1.\end{equation}
\noindent
For all $\lambda\in\bn$ consider the $n\times(n+2)$-matrix with integer entries
$$A_{\lambda}=\left(
\begin{array}{cccccc}
\lambda a_1&b_1&c&0&\cdots&0\\
\vdots&\vdots&0&\ddots&&\vdots\\
\vdots&\vdots&\vdots&&\ddots&0\\
\lambda a_n&b_n&0&\cdots&0&c
\end{array}
\right),
$$
\noindent
and the $n\times(n+1)$-submatrix
$$A'=\left(
\begin{array}{ccccc}
b_1&c&0&\cdots&0\\
\vdots&0&\ddots&&\vdots\\
\vdots&\vdots&&\ddots&0\\
b_n&0&\cdots&0&c
\end{array}
\right).
$$
\noindent
For all $\lambda\in\bn$ let $g_{\lambda}$ be the greatest common divisor of the $n$-minors of $A_{\lambda}$ and let $g'$ be the greatest common divisor of the $n$-minors of $A'$. Moreover, set
$$\bar g=\gcd_{1\leq i<j\leq n}(c^{n-1}a_i, c^{n-2}(a_ib_j-a_jb_i)).$$
\noindent
Then
$$g_{\lambda}=\gcd(g',\lambda\bar g).$$
\noindent
If we apply the
 solvability criterion for diophantine systems given in \cite{He} to  $A'{\bf x}=\lambda{\bf a}$, we conclude that $\lambda{\bf a}\in\bz\tilde T_2$ if and only if $g'=g_{\lambda}$, i.e., if and only if $g'$ divides $\lambda \bar g$, which is the case if and only if $\lambda$ is a multiple of
\begin{equation}\label{alpha1} \alpha=\frac{g'}{\gcd(g', \bar g)}.\end{equation}
\noindent
Note that
\begin{equation}\label{g'}g'=\gcd_{1\leq i\leq n}(c^n, c^{n-1}b_i)=c^{n-1}\gcd_{1\leq i\leq n}(c, b_i)=c^{n-2}\gcd_{1\leq i\leq n}(c^2, cb_i),\end{equation}
\noindent and
$$\bar g=c^{n-2}\gcd_{1\leq i<j\leq n}(ca_i, a_ib_j-a_jb_i),$$
\noindent
so that
\begin{eqnarray}\label{gg}
\gcd(g',\bar g)&=&c^{n-2}\gcd_{1\leq i<j\leq n}(c^2, ca_i, cb_i, a_ib_j-a_jb_i)\nonumber\\
&=&c^{n-2}\gcd_{1\leq i<j\leq n}(c\gcd(c, a_i, b_i), a_ib_j-a_jb_i)\nonumber\\
&=&c^{n-2}\gcd_{1\leq i<j\leq n}(c, a_ib_j-a_jb_i),
\end{eqnarray}
\noindent
where the last equality follows from (\ref{gcd}).
Finally, from (\ref{alpha1}), (\ref{g'}) and (\ref{gg}) we deduce that
\begin{equation}\label{alpha2} \alpha=\frac{c\gcd_{1\leq i\leq n}(c, b_i)}{\gcd_{1\leq i\leq n}(c, a_ib_j-a_jb_i)},\end{equation}
\noindent
which is the sought closed formula. 
\par\smallskip\noindent
We are now ready to prove the main result of this section, which completes Corollary 2 in \cite{BMT1}.
\begin{proposition}\label{gT} Suppose that neither of supp\,${\bf a}$ and supp\,${\bf b}$ is contained in the other. Let $c'=\gcd_{1\leq i<j\leq n}(c, a_ib_j-a_jb_i)$, and set
$$g(T)=\frac{c'}{\gcd_{1\leq i\leq n}(c', a_i)\gcd_{1\leq i\leq n}(c', b_i)}.$$
\noindent
Then $T$ is completely $p$-glued
\begin{list}{}{}
\item{(a)} for exactly one prime $p$, if and only if $g(T)$ is a positive power of $p$;
\item{(b)} for  every prime $p$, if and only if $g(T)$ is equal to 1;
\item{(c)} for no prime $p$, if and only if $g(T)$ has two distinct prime divisors.
\end{list}
Case (b) occurs if and only if $T$ is a complete intersection.
\end{proposition}
\demo We preliminarily note that the last claim is a consequence of Corollary \ref{corollary6bis}. Let $p$ be any prime. By virtue of Theorem \ref{main} and Proposition \ref{proposition2} b), $T$ is completely $p$-glued if and only if it is the $p$-gluing of $\tilde T_1$ and $\tilde T_2$. This occurs if and only if $p^k\alpha{\bf a}\in\bn\tilde T_2$, which is equivalent to
\begin{equation}\label{pka} p^k\alpha{\bf a}\in\bn\{c{\bf e}_1,\dots, c{\bf e}_n\},\end{equation}
\noindent
for some nonnegative integer $k$, where $\alpha$ is the positive integer defined in (\ref{alpha}) and given by the formula (\ref{alpha2}). Now (\ref{pka}) is true if and only if $p^k$ is a multiple of
$$\omega=\frac{c}{\gcd_{1\leq i\leq n}(c, \alpha a_i)}=\frac{c}{\gcd(c, \alpha\gcd_{1\leq i\leq n}(a_i))}.$$
\noindent
Thus (\ref{pka}) is true for exactly one prime $p$, for every prime $p$, or for no prime $p$ if $\omega$ is a positive power of $p$, is equal to 1, or has two distinct prime divisors, respectively. In view of (\ref{alpha2}) we have
\begin{eqnarray*}
\omega&=&\frac{c}{\gcd(c, \frac{c}{c'}\gcd_{1\leq i\leq n}(c, b_i)\gcd_{1\leq i\leq n}(a_i))}\\
&=&\frac{c}{\frac{c}{c'}\gcd(c',\gcd_{1\leq i\leq n}(c, b_i)\gcd_{1\leq i\leq n}(a_i))}\\
&=&\frac{c'}{\gcd(c',\gcd_{1\leq i\leq n}(c', b_i)\gcd_{1\leq i\leq n}(a_i))}\\
&=&\frac{c'}{\gcd(c',\gcd_{1\leq i\leq n}(c', b_i)\gcd_{1\leq i\leq n}(c',a_i))}\\
&=&\frac{c'}{\gcd_{1\leq i\leq n}(c', b_i)\gcd_{1\leq i\leq n}(c',a_i)},
\end{eqnarray*}
\noindent
where the last equality follows from the fact that, as a consequence of (\ref{gcd}), $\gcd(c',a_1,\dots, a_n, b_1, \dots, b_n)=1$. This shows that $\omega=g(T)$ and completes the proof.
\begin{remark}{\rm Since $g(T)$ is a divisor of $c$, case (c) of Proposition \ref{gT} cannot hold whenever $c$ is a prime power. From this we deduce the well-known fact that $T$ is always completely $p$-glued if $c=p^r$ for some nonnegative integer $r$ (this statement holds for simplicial toric varieties of any codimension: for a proof see, e.g., \cite{B4}, Proposition 2). In \cite{BL} we considered the case where $c=p$. There, in Theorem 2.1, we gave a different characterization of case (a)  in form of four arithmetic conditions.} \end{remark}
\section{Applications to set-theoretic complete intersections}
We now use the results on $p$-gluings proven in Section 3 for giving a complete characterization of  simplicial toric varieties of codimension 2 with respect to the property of being a binomial set-theoretic complete intersection. The (affine or projective) toric variety $V_T$ attached to the set $T$ considered in the previous sections admits the following parametrization:
$$x_1=u_1^c,\quad \dots,\quad x_n=u_n^c,\quad y_1=u_1^{a_1}\cdots u_n^{a_n},\quad y_2=u_1^{b_1}\cdots u_n^{b_n}.$$
\noindent
As an immediate consequence of Proposition \ref{p} and Corollary \ref{corollary_added} we have the following.
\begin{corollary}\label{corollary7} A simplicial toric variety of codimension 2 is a binomial set-theoretic complete intersection in characteristic $p$ either
\begin{list}{}{}
\item{(a)} for exactly one prime $p$, or
\item{(b)} for  every prime $p$, or
\item{(c)} for no prime $p$.
\end{list} 
\end{corollary}
\noindent
\begin{example}{\rm 
Suppose that $T$ is the $p$-gluing of $\tilde T_1$ and $\tilde T_2$, with $\bz\tilde T_1\cap\bz\tilde T_2=\bz\alpha{\bf a}$ and $p^k\alpha{\bf a}=\beta{\bf b}+\sum_{i=1}^n\beta_ic{\bf e}_i$, where $k,\beta, \beta_1\,\dots \beta_n\in\bn$. Then the construction performed in the proof of \cite{BMT2}, Theorem 5, shows that $V_T$ is a set-theoretic complete intersection defined by the vanishing of the following two binomials:
$$F_1=y_1^{p^k\alpha}-y_2^{\beta}x_1^{\beta_1}\cdots x_n^{\beta_n},\qquad  F_2=y_2^{\gamma}-x_1^{\delta_1}\cdots x_n^{\delta_n}.$$
Here we have set $\gamma=\frac{c}{\gcd_{1\leq j\leq n}(c,b_j)}$, and $\delta_i=\frac{b_i}{\gcd_{1\leq j\leq n}(c,b_j)}$ for all $i=1,\dots, n$.}
\end{example}
\begin{remark}\label{Lyubeznik}{\rm According to \cite{BL}, Theorem 2.1, if $c=p$, then $V_T$ is a set-theoretic complete intersection in some positive characteristic if and only if it is a binomial set-theoretic complete intersection in the same characteristic. } 
\end{remark} 
From Corollary \ref{corollary1} we know that, if the dimension $n$ of $V$ is at most 2, case (b) of Corollary \ref{corollary7} always occurs. On the contrary, for every $n\geq 3$, there are
 simplicial toric varieties of codimension 2 for each of the cases (a)--(c). On the other hand, Proposition \ref{full} tells us that case (b) is fulfilled whenever the parametrization of $V$ is full. We give examples of the cases (a) and (c) for $n=3$.  
\begin{example}{\rm Let $p$ be a prime and let $V_T\subset K^5$ be the affine toric variety attached to the set 
$$T=\{(p,0,a), (0,p,b), (p,0,0), (0,p,0), (0,0,p) \},$$
\noindent
where $a$ and $b$ are positive integers not divisible by $p$. 
Let
$$\tilde T_1=\{(p,0,a)\}, \qquad\tilde T_2=\{(0,p,b), (p,0,0), (0,p,0), (0,0,p)\}.$$
\noindent
 Then
$$\bz(p,0,a)=\bz \tilde T_1\cap\bz\tilde T_2.$$
\noindent
In fact, if $\lambda$ and $\mu$ are integers such that $\lambda b +\mu p=a$, then
$$(p,0,a)=\lambda(0,p,b)+(p,0,0)-\lambda(0,p,0)+\mu(0,0,p).$$
\noindent
On the other hand, for every positive integer $\alpha$, we have that $\alpha(p,0,a)\in\bn\tilde T_2$ if and only if $\alpha(p,0,a)\in\bn\{(p,0,0), (0,0,p)\}$, which occurs if and only if $p$ divides $\alpha$. This shows that $T$ is the $p$-gluing of $\tilde T_1$ and $\tilde T_2$, but not the $q$-gluing of  $\tilde T_1$ and $\tilde T_2$ for any prime $q\ne p$. In view of Theorem \ref{main}, this implies that  $T$ is completely $p$-glued, but not completely $q$-glued for any prime $q\ne p$. Hence $p$ is the only positive characteristic in which $V_T$ is a binomial set-theoretic complete intersection; in particular,  $V_T$ is not a complete intersection, hence, in view of Proposition \ref{Barile}, it is not a binomial set-theoretic complete intersection in characteristic zero. In characteristic $p$, the variety $V_T$ is defined by the vanishing of the binomials
$$F_1=y_1^p-x_1^px_3^a,\qquad\qquad F_2=y_2^p-x_2^px_3^b.$$
\noindent 
According to Remark \ref{Lyubeznik}, $V_T$ is not a set-theoretic complete intersection in any characteristic other than $p$. 
}
\end{example}
\begin{example}{\rm 
Let $p,q$ be different primes and let $V\subset K^5$ be the affine toric variety attached to the set 
$$T=\{(pq,0,a), (0,pq,b), (pq,0,0), (0,pq,0), (0,0,pq) \},$$
\noindent
where $a$ and $b$ are positive integers not divisible by $p$ nor $q$.
Let
$$\tilde T_1=\{(pq,0,a)\}, \qquad\tilde T_2=\{(0,pq,b), (pq,0,0), (0,pq,0), (0,0,pq)\}.$$
\noindent
 Then
$$\bz(pq,0,a)=\bz \tilde T_1\cap\bz\tilde T_2.$$
\noindent
In fact, if $\lambda$ and $\mu$ are integers such that $\lambda b +\mu pq=a$, then
$$(pq,0,a)=\lambda(0,pq,b)+(pq,0,0)-\lambda(0,pq,0)+\mu(0,0,pq).$$
\noindent
On the other hand, for every positive integer $\alpha$, we have that $\alpha(pq,0,a)\in\bn\tilde T_2$ if and only if $\alpha(pq,0,a)\in\bn\{(pq,0,0), (0,0,pq)\}$, which occurs if and only if $pq$ divides $\alpha$. This shows that $T$ is the not the $p$-gluing of $\tilde T_1$ and $\tilde T_2$ for any prime $p$. In view of Theorem \ref{main}, this implies that  $T$ is  not completely $p$-glued for any prime $p$. Hence  $V_T$ is not a binomial set-theoretic complete intersection in any characteristic.\newline
The variety $V_T$ belongs to the class considered in \cite{B4}, Section 3 (see, in particular, Example 5). There it was shown that $V_T$ is in fact not a set-theoretic complete intersection in any characteristic. 
}
\end{example}
\noindent
Proposition \ref{proposition1} and Corollary \ref{corollary_added} imply that, if $T$ is completely $p$-glued for some prime $p$, then $T$ is the $q$-gluing of $\tilde T_1$ and $\tilde T_2$ (or of $\tilde{\tilde{T_1}}$ and $\tilde{\tilde{T_2}}$) either for every prime $q$, or for no prime $q$ or only for $q=p$. The situation is different for higher codimensions. The next example, which refers to codimension 3, shows that the set 
$$T=\{{\bf a}, {\bf b}, {\bf c}, {\bf e}_1, \dots, {\bf e}_n\}$$ 
can admit two different partitions 
$$\tilde T_1=\{{\bf a}\}, \qquad\tilde T_2=\{{\bf b}, {\bf c}, {\bf e}_1, \dots, {\bf e}_n\},$$
$$\tilde{\tilde T_1}=\{{\bf b}\}, \qquad\tilde{\tilde T_2}=\{{\bf a}, {\bf c}, {\bf e}_1, \dots, {\bf e}_n\},$$
\noindent
with respect to which it is only $p$-glued and only $q$-glued, respectively, where $p$ and $q$ are two distinct primes.
\begin {example}{\rm Let  $V_T\subset K^6$ be the affine toric variety attached to the set 
$$T=\{(3, 0, 3), (0,4,2), (3,2,1), (6,0,0), (0,6,0), (0,0,6)\},$$
\noindent
and let
$$\tilde T_1=\{(3,0,3)\}, \qquad\tilde T_2=\{(0,4,2), (3,2,1), (6,0,0), (0,6,0), (0,0,6)\}.$$
\noindent
Then $\tilde T_2$ is completely $p$-glued for all primes $p$ by Proposition \ref{proposition1}. Moreover, 
$$(3,0,3)=(6,0,0)+(0,6,0)+(0,0,6)-(3,2,1)-(0,4,2),$$
whence we deduce that
$$\bz(3,0,3)=\bz \tilde T_1\cap\bz\tilde T_2.$$
\noindent
For every positive integer $\alpha$, we have that $\alpha(3,0,3)\in\bn\tilde T_2$ if and only if $\alpha(3,0,3)\in\bn\{(6,0,0), (0,0,6)\}$, which occurs if and only if $2$ divides $\alpha$. This shows that $T$ is the 2-gluing of $\tilde T_1$ and $\tilde T_2$, so that $T$ is completely 2-glued; and it is not the $p$-gluing  of $\tilde T_1$ and $\tilde T_2$ for any prime $p\ne2$. However, if
$$\tilde{\tilde T_1}=\{(0,4,2)\}, \qquad\tilde{\tilde T_2}=\{(3,0,3), (3,2,1), (6,0,0), (0,6,0), (0,0,6)\},$$
\noindent
then $T$ is the 3-gluing of $\tilde{\tilde T_1}$ and $\tilde{\tilde T_2}$, so that $T$ is also completely 3-glued; but it is not the $p$-gluing  of $\tilde T_1$ and $\tilde T_2$ for any prime $p\ne3$. \newline
Finally, it turns out that, for all primes $p$, $T$ is the $p$-gluing of 
$$\bar T_1=\{(3,2,1)\}, \qquad\bar T_2=\{(3,0,3), (0,4,2), (6,0,0), (0,6,0), (0,0,6)\},$$
\noindent
where $\bar T_2$ is completely $p$-glued. Hence $T$ is completely $p$-glued for all primes $p$. \newline
We do not know of any affine toric variety of codimension at least 3 whose attached set is completely $p$-glued for two different primes $p$, but not for all primes $p$. The above example demonstrates that proving that such varieties do not exist in codimension greater than 2 would require arguments which differ substantially from those we have applied to the case of codimension 2. 
}
\end{example}


\begin{thebibliography}{BMT2} 
\bibitem{B0} M.~Barile, Arithmetical ranks of ideals associated to symmetric and alternating matrices,
{\it J. Algebra} {\bf 176} (1995), 59--82.
\bibitem{B1} M.~Barile, A note on Veronese varieties. {\it Rend. Circ. Mat. Palermo} {\bf 54} (2005), 359--366.
\bibitem{B2} M.~Barile, On toric varieties of high arithmetical rank. {\it Yokohama Math. J.} {\bf 52} (2006), 125--130. 
\bibitem{B3} M.~Barile, On a special class of simplicial toric varieties. {\it J. Algebra} {\bf 308} (2007), 368--382.
\bibitem{B4} M.~Barile, On simplicial toric varieties of codimension 2. Preprint (2007). math.AC/0705.4389. To appear in: {\it Rend. Istit. Mat. Univ. Trieste}. 
\bibitem{BL} M.~Barile, G.~Lyubeznik, Set-theoretic complete intersections in characteristic $p$. {\it Proc. Amer. Math. Soc.} {\bf 133} (2005), 3199--3209.
\bibitem{BMT1} M.~Barile, M.~Morales, A.~Thoma, On simplicial toric varieties which are set-theoretic complete intersections. {\it J. Algebra} {\bf 226} (2000), 880--892.
\bibitem{BMT2} M.~Barile, M.~Morales, A.~Thoma, Set-theoretic complete intersections on binomials. {\it Proc. Amer. Math. Soc.} {\bf 130} (2002), 1893--1903.
\bibitem{Br1} H.~Bresinsky, Monomial space curves in ${\bf A}\sp{3}$ as set-theoretic complete intersections. {\it Proc. Amer. Math. Soc.} {\bf 75} (1979),  23--24.
\bibitem{Br2} H.~Bresinsky, Monomial Gorenstein curves in ${\bf A}\sp{4}$ as set-theoretic complete intersections. {\it Manuscripta Math.} {\bf 27} (1979), 353--358.
\bibitem{E} K.~Eto, An example of set-theoretic complete intersection lattice
ideal, {\it Tokyo Math. J.} {\bf 29} (2006), 319–-324.
\bibitem{FMS} K.~Fischer, W.~Morris, J.~Shapiro, Affine semigroups that are complete intersections. {\it Proc. Amer. Math. Soc.} {\bf 125} (1997), 3137--3145.
\bibitem{He} J.~Heger, Denkschriften. {\it Kais. Akad. Wissensch. Mathem. Naturw. Klasse} {\bf 14} (1858), II.
\bibitem{K} A.~Katsabekis, Projection of cones and the arithmetical rank of toric
varieties. {\it J. Pure Appl. Algebra} {\bf 199} (2005), 133–-147.
\bibitem{M} T.T.~Moh, Set-theoretic complete intersections. {\it Proc. Amer. Math. Soc.} {\bf 94} (1985), 217--220.
\bibitem{Mo} M.~Morales, Equations of monomial varieties in codimension two.  {\it J. Algebra} {\bf 175} (1995), 1082--1095.
\bibitem{RV1} L.~Robbiano, G.~Valla, On set-theoretic complete intersections in the projective space. {\it Rend. Sem. Mat. Fis. Milano} {\bf 53} (1983), 333--346. 
\bibitem{RV2} L.~Robbiano, G.~Valla, Some curves in ${\bf P}\sp{3}$ are set-theoretic complete intersections. Algebraic geometry---open problems (Ravello, 1982), 391--399, Lecture Notes in Math., {\bf 997}, Springer, Berlin-New York, 1983. 
\bibitem{R} J.C.~Rosales, On presentations of subsemigroups of ${\bn}^n$. {\it Semigroup Forum} {\bf 55} (1997), 152--159.
\bibitem{S} M.~\c{S}ahin, Producing complete intersection monomial curves in ${\bf P}^n$. Preprint (2006). math.AG/0610316.
\bibitem{T} A.~Thoma, On the set-theoretic complete intersection problem for
monomial curves in ${\bf A}^n$ and ${\bf P}^n$, {\it J. Pure Appl. Algebra} {\bf 104} (1995), 333–-344.
\bibitem{V}  G.~Valla, On determinantal ideals which are set-theoretic complete intersections. {\it Compositio Math.} {\bf 42} (1980/81), 3--11.
\end{thebibliography}
\end{document}